\documentclass[a4paper]{amsart}

\usepackage[T1]{fontenc}
\usepackage{lmodern}
\usepackage[pdfusetitle]{hyperref}

\def\Z{\mathbb Z}
\def\Q{\mathbb Q}
\def\R{\mathbb R}
\def\C{\mathbb C}
\def\bP{\mathbb P}
\def\CP{\mathbb{CP}}
\def\cO{\mathcal O}
\newcommand{\PP}{\mbox{$\mathit P\hspace{-.8pt}P$}}
\DeclareMathOperator{\lcm}{lcm}
\DeclareMathOperator{\codim}{codim}

\theoremstyle{plain}
\newtheorem{theorem}{Theorem}[section]
\newtheorem{proposition}[theorem]{Proposition}
\newtheorem{lemma}[theorem]{Lemma}

\theoremstyle{definition}

\newtheorem{remark}[theorem]{Remark}

\newenvironment{acknowledgements}%
  {\begin{trivlist}\item[]{\em Acknowledgements.}\hskip0.5em}
  {\end{trivlist}}

\numberwithin{equation}{section}

\begin{document}

\title[The equivariant cohomology of weighted projective space]%
  {The equivariant cohomology of\\weighted projective space}
\author{Anthony Bahri, Matthias Franz and Nigel Ray}

\begin{abstract}
  We describe the integral equivariant cohomology ring of a weighted
  projective space in terms of piecewise polynomials, and thence by
  generators and relations. We deduce that the ring is a perfect
  invariant, and prove a Chern class formula for weighted projective
  bundles.
\end{abstract}


\maketitle

%
%
%
%
%
%
%
%
%

\section{Introduction}\label{sec:1}

Let $\chi=(\chi_0,\ldots,\chi_n)$ be a vector of positive natural
numbers.  The associated \emph{weighted projective space} is the
quotient
\begin{equation}\label{definition-weighted-projective-space}
  \bP(\chi)=S^{2n+1}/S^1\langle\chi_0,\ldots,\chi_n\rangle,
\end{equation}
where the numbers $\chi_i$ indicate the weights with which $S^1$ acts
on the unit sphere $S^{2n+1}\subset\C^{n+1}$, by
\begin{equation}
  g\cdot(x_0,\ldots,x_n)=(g^{\chi_0}x_0,\ldots,g^{\chi_n}x_n).
\end{equation}
So $\bP(1,\dots,1)$ is the standard projective space $\CP^n$. 
Note that $\bP(\chi)$ is equipped with an action of the
$n$-dimensional torus
\begin{equation}\label{wtdtorus}
  T=(S^1)^{n+1}/S^1\langle\chi_0,\ldots,\chi_n\rangle,
\end{equation}
where the quotient is defined by analogy with
(\ref{definition-weighted-projective-space}).  We give an explicit
example of such an action in (\ref{act1234}), for $n=3$.  Different
weight vectors may also give equivalent weighted projective spaces; we
will elaborate on this aspect in Section~\ref{recover-weights}.

Kawasaki~\cite{Kawasaki:73} has computed the ordinary cohomology ring
of~$\bP(\chi)$ with integer coefficients. Additively, the cohomology
is isomorphic to that of $\CP^n$, but the multiplication is different.
More precisely, if $c_1$ is a generator of the group $H^2(\bP(\chi))$,
then $H^*(\bP(\chi))$ is generated as a ring by elements $c_m$, where
$1\leq m\leq n$ and
\begin{equation}\label{divisibility-powers} 
  c_1^m=\frac{\lcm(\chi_0,\dots,\chi_n)^m}
  {\lcm\{\,\prod_{i\in I}\chi_i:|I|=m\,\}}\,c_m
\end{equation}
in $H^{2m}(\bP(\chi))$. The multiplication is induced accordingly.

In this note we study $H_T^*(\bP(\chi))$, the $T$-equivariant
cohomology of~$\bP(\chi)$ with integer coefficients; it is defined as
the ordinary cohomology $H^*(\bP(\chi)_T)$ of the Borel construction
$\bP(\chi)_T=ET\times_T\bP(\chi)$.  Our main result,
Theorem~\ref{wproj-generators}, describes $H_T^*(\bP(\chi))$ in terms
of generators and relations.

We give two applications. Firstly, we show how to recover the weight
vector $\chi$ from $H_T^*(\bP(\chi))$, thereby establishing that
different weighted projective spaces have different integral
equivariant cohomology rings.  Secondly, we consider weighted
projective bundles.  For any direct sum $D=L_1\oplus\dots\oplus L_n$
of complex line bundles over a base space $X$, the cohomology ring
$H^*(\bP(D))$ of the projectivisation is a module over $H^*(X)$.  Its
algebra structure is determined by the single relation
\begin{equation}\label{chern-projective-prod}
  \prod_{i=0}^n\bigl(\xi+c_1(L_i)\bigr) = 0,
\end{equation}
where $c_1(L_i)\in H^2(X)$ and $-\xi\in H^2(\bP(D))$ denote the Chern
classes of $L_i$ and of the canonical complex line bundle
respectively. Al~Amrani~\cite{AlAmrani:97} has stated a generalisation
of (\ref{chern-projective-prod}) to weighted projective bundles and
proved it in a special case. Theorem~\ref{weighted-chern} establishes
his relation in general.

We consider our calculations as lying in the realm of toric topology,
and will elaborate on this theme in a subsequent document.  Readers
who require background information on equivariant topology may consult
\cite{AlldayPuppe:93}, or the survey articles in \cite{May:96}.

\begin{acknowledgements}
The first author thanks Rider University for the award of Research
Leave. All three authors are grateful to the Manchester Institute of
Mathematical Sciences (MIMS) for ongoing support that has helped to
sustain their collaboration.
\end{acknowledgements}

%
%
%
%
%
%
%
%
%

\section{From equivariant cohomology to piecewise polynomials}
\label{to-piecewise}

By a \emph{ring} we always mean a graded commutative ring with unit
element. All rings we consider happen to be concentrated in even
degrees, so that they are commutative in the ordinary sense.

\begin{remark}
  There are several ways to describe the divisibility of the 
  powers $c_1^m$ in $H^*(\bP(\chi))$. Kawasaki looks at the
  $p$-contents of the weights for each prime $p$ separately. If
  $q_0(p)$,~\ldots,~$q_n(p)$ are their $p$-contents, in increasing
  order, then
 \begin{equation}\label{divisibility-kawasaki} 
    c_1^m = \prod_p 
            \frac{q_n(p)^m}{q_n(p)\cdots q_{n-m+1}(p)}\:c_m. 
  \end{equation}
  Kawasaki also considers sets~$I$ of size~$m+1$ instead of~$m$ and
  writes
  \begin{equation}\label{divisibility-alamrani}
    c_1^m = \frac{\lcm(\chi_0,\dots,\chi_n)^m} 
  {\lcm\{\,\gcd\{\chi_i: i\in I\}^{-1} 
  \prod_{i\in I}\chi_i:|I|=m+1\,\}} 
    \,c_m 
  \end{equation}
  \cite[p.248]{Kawasaki:73},  
  as does Al~Amrani~\cite[Sec.~I.5]{AlAmrani:97}.
  Taking products of $m+1$~weights and
  dividing by their greatest common divisor,
  as in the denominator of~(\ref{divisibility-alamrani}),
  removes the smallest $p$-content for each prime~$p$. Computing
  the least common multiple over all such terms then gives the product
  of the $m$~largest $p$-contents of all weights, in accordance
  with the denominators of 
  (\ref{divisibility-powers})~and~(\ref{divisibility-kawasaki}).
\end{remark}
Now let $\iota\colon\bP(\chi)\to \bP(\chi)_T$ be the inclusion of a fibre
into the Borel construction.
\begin{lemma}\label{HT-free-generators}
  As an $H^*(BT)$-module, $H_T^*(\bP(\chi))$ is free of rank $n+1$: as
  a ring, it is generated by the image of $H^2(BT)$
  in~$H_T^2(\bP(\chi))$, together with any subgroup
  $A^*<H_T^*(\bP(\chi))$ that surjects onto $H^{>0}(\bP(\chi))$
  under $\iota^*$.
\end{lemma}
\begin{proof}
  By Kawasaki, $H^*(\bP(\chi))$ is free over~$\Z$ and concentrated in
  even degrees.  So the Serre spectral sequence of the
  fibration $\bP(\chi)\to\bP(\chi)_T\to BT$ degenerates at the
  $E_2$ level, and $H_T^*(\bP(\chi))$ is isomorphic as
  $H^*(BT)$-modules to $H^*(\bP(\chi))\otimes H^*(BT)$ by the
  Leray--Hirsch theorem. The isomorphism is induced by any additive
  section to $\iota^*$, which we may assume takes values in $A^*$.
  This proves the claim.
\end{proof}

The equivariant cohomology of $\CP^n$ is well-known, and may be
described conveniently in the context of toric varieties. Indeed,
$\bP(\chi)$ is an $n$-dimensional projective toric variety for every
$\chi$, and may be constructed from any complete simplicial fan~%
$\Sigma$ whose rays $v_0$,~\ldots,~$v_n\in N=\Z^n$ span $N$, and
satisfy the relation
  \begin{equation}
    \label{relation-fan-generators}
    \chi_0 v_0 + \cdots + \chi_n v_n = 0.
  \end{equation}
Every such $\Sigma$ is the normal fan of a lattice $n$-simplex in
$N\otimes_\Z\R$, and is polytopal; we refer to
Fulton~\cite[\S2.2]{Fulton:93}, or the nice overview
in~\cite[\S4.1]{Kasprzyk:06}, for a full discussion.

In particular, $H_T^*(\CP^n)$ is isomorphic to the integral 
Stanley--Reisner algebra 
\begin{equation}\label{stanreisbdrysplx}
  \Z[\Sigma] = \Z[a_0,\ldots,a_n]/(a_0\cdots a_n)
\end{equation}
of the appropriate $\Sigma$, where each generator~$a_i$ corresponds to
$v_i$, and has degree~$2$. In other words, the only relation amongst
the generators is
\begin{equation}\label{relation-product-all}
  \prod_{i=0}^n a_i = 0.
\end{equation}
The situation for singular varieties is less straightforward, and 
the $\bP(\chi)$ offer a natural family of test cases. Our aim is to 
generalise (\ref{stanreisbdrysplx}), and express $H_T^*(\bP(\chi))$ 
in terms of generators and relations for arbitrary $\chi$. We use 
the language of piecewise polynomials, to which we now turn.

A function~$f\colon N\to\Z$ is called \emph{piecewise polynomial} on
$\Sigma$ if it coincides with some globally defined polynomial
$g\in\Z[N]$ on each cone $\sigma$. Such functions are closed under
pointwise addition and multiplication, and form an algebra
$\PP\hspace{1pt}[\Sigma]$ over the ring of global polynomials $\Z[N]$.
We grade both $\PP[\Sigma]$ and $\Z[N]$ by twice the degree of
homogeneous elements, and use the Chern classes of the $n$ canonical
line bundles to identify $\Z[N]$ with $H^*(BT)$.

\begin{proposition}\label{iso-HT-piecewise-poly}
Let $\Sigma$ be a polytopal fan in $N$, and $X_\Sigma$ 
the associated compact projective toric variety: if 
$H^*(X_\Sigma)$ is concentrated in even degrees, then 
$H_T^*(X_\Sigma)$ is isomorphic to 
$\PP[\Sigma]$ as an algebra over $H^*(BT)$. 
\end{proposition}
\begin{proof}
  Set $X=X_\Sigma$, and denote the orbit space by $X/T$. The latter
  may be identified with a convex polytope $P_\Sigma$, whose normal
  fan is $\Sigma$. Following Goresky and MacPherson
  \cite{GoreskyMacPherson:87}, we deduce that $X$ is homeomorphic to a
  quotient space of~$T\times P_\Sigma$, and may therefore be expressed
  as a finite $T$-CW complex \cite{May:96} with connected isotropy
  groups.  As in the proof of Lemma~\ref{HT-free-generators}, the
  Serre spectral sequence for the fibration~$X\to X_T\to BT$
  degenerates at the $E_2$~level because $H^*(X)$ is concentrated in
  even degrees.  Hence, by a result of
  Franz--Puppe~\cite{FranzPuppe:exact}, the Chang--Skjelbred sequence
  \begin{equation}\label{chang-skjelbred}
    0 \longrightarrow H^*_T(X)
      \stackrel{j^*}\longrightarrow H^*_T(X^T)
      \stackrel{\delta}\longrightarrow H^{*+1}_T(X_1, X^T)
  \end{equation}
  is exact for integral coefficients.  Here $X^T$ denotes the
  $T$-fixed points, $X_1$ the union of~$X^T$ and all $1$-dimensional
  orbits, $j$ the inclusion~$X^T\to X$ and $\delta$ the differential
  of the long exact cohomology sequence for the pair~$(X_1,X^T)$.  We
  may identify the kernel of $\delta$ with the algebra $\PP[\Sigma]$,
  as follows.

  Write $\cO_\sigma$ for the orbit under the complexification~$T_\C$
  of~$T$ corresponding to the cone $\sigma\in\Sigma$, and $\Z[\sigma]$
  for the polynomials with integer coefficients on the linear hull
  of~$\sigma$. Any such polynomial is uniquely determined by its
  restriction to~$\sigma$.

  When $\sigma$ is $n$-dimensional, we have that
  \begin{equation}
    H_T^*(\cO_\sigma)=H^*(BT)=\Z[\sigma].
  \end{equation}
  For $(n-1)$-dimensional $\tau$, we denote the isotropy group
  of~$\cO_\tau$ by~$T_\tau$. Then the action of the circle~$T/T_\tau$
  on the closure~$\bar\cO_\tau$ is isomorphic to the standard action
  of~$S^1$ on~$\CP^1$, whose fixed points we write as
  $0$~and~$\infty$. We obtain
\begin{equation}
    H_T^*(\bar\cO_\tau,\partial\cO_\tau)
    = H^*(BT_\tau)\otimes H_{T/T_\tau}^*(\bar\cO_\tau,\partial\cO_\tau)
\end{equation}
 and 
 \begin{equation}
    H_{T/T_\tau}^*(\bar\cO_\tau,\partial\cO_\tau)
    \cong H_{S^1}^*(\CP^1,\{0,\infty\})
    \cong \Z[+1];
\end{equation}
  hence
\begin{equation}
    H_T^{*+1}(\bar\cO_\tau,\partial\cO_\tau)
    \cong \Z[\tau].
  \end{equation}
  Moreover, whenever $\tau\subset\sigma$ is a facet, the differential
  \begin{equation}
    H_T^*(\cO_\sigma)\longrightarrow
    H_T^{*+1}(\bar\cO_\tau,\partial\cO_\tau)
  \end{equation}
  is the canonical restriction~$\Z[\sigma]\to\Z[\tau]$, multiplied
  by~$\pm1$ according to the orientation of the
  interval~$\cO_\tau/T\approx(0,\infty)$.

  It follows that the differential~$\delta$ of (\ref{chang-skjelbred})
  is a signed sum of restrictions
  \begin{equation}
    \delta\colon \bigoplus_{\sigma\in\Sigma^n}\Z[\sigma]
    \longrightarrow 
    \bigoplus_{\tau\in\Sigma^{n-1}}\Z[\tau]
  \end{equation}
  of sums of polynomial algebras, whose component into~$\Z[\tau]$ is
  the difference of the restrictions of the polynomials on the two
  $n$-dimensional cones having $\tau$ as their common facet. So the
  kernel consists of those collections of polynomials on
  $n$-dimensional cones which may be glued along their common facets.
  This corresponds to the requirement that the polynomials agree on
  \emph{any} intersection~$\tau=\sigma\cap\sigma'$, because
  $\sigma$~and~$\sigma'$ are connected by a sequence of
  $n$-dimensional cones, each of which contains $\tau$ and shares a
  facet with the next.  (In other words, $\Sigma$ is a
  \emph{hereditary fan} \cite{BilleraRose:92}.)  Thus $\delta$ has
  kernel $\PP[\sigma]$, as required.

  Finally, let $f$ be the piecewise polynomial corresponding to a
  class $\alpha\in H_T^*(X)$. By construction, the polynomial
  coinciding with $f$ on an $n$-dimensional cone $\sigma$
  is the image of~$\alpha$ under the restriction
  \begin{equation}
  H_T^*(X)\to H_T^*(\cO_\sigma)=H^*(BT)=\Z[N].
  \end{equation}
  Since the composition $H^*(BT)\to H_T^*(X) 
  \to H_T^*(\cO_\sigma)=H^*(BT)$
  is the identity, it follows that the map~$H^*(BT)\to H_T^*(X)$
  corresponds to the inclusion of the algebra of global polynomials.
\end{proof}

\begin{remark}
  The integral equivariant cohomology of any smooth, not necessarily
  compact toric variety~$X_\Sigma$ is given by the Stanley--Reisner
  algebra $\Z[\Sigma]$ \cite{BifetDeConciniProcesi:90},
  \cite{DavisJanuszkiewicz:91}; or equivalently, by $\PP[\Sigma]$
  \cite{Brion:96}.  A canonical isomorphism between the two is defined
  by assigning the Courant function $a_\rho$ of the ray~$\rho$ to the
  Stanley--Reisner generator corresponding to~$\rho$. The function
  $a_\rho$ is piecewise linear on~$\Sigma$, and assumes the value~$1$
  on the generator of~$\rho$ and $0$ on all other rays. It is
  well-defined because the smoothness of~$X_\Sigma$ implies that the
  rays of any cone may be completed to a basis of the lattice~$N$.
  Brion was probably the first to note the relationship between
  piecewise polynomials and the Chang--Skjelbred sequence.

  Similarly, when $\Sigma$ is simplicial the rational equivariant
  cohomology~$H_T^*(X_\Sigma;\Q)$ is the Stanley--Reisner
  algebra~$\Q[\Sigma]$; or equivalently, $\PP[\Sigma]\otimes\Q$
  \cite[p.107]{Fulton:93}.  In particular, there is an isomorphism
\begin{equation}\label{wpsratcoho} 
  H^*_T(\bP(\chi);\Q)\cong\Q[a_0,\dots,a_n]/(a_0\dots a_n).  
\end{equation}  

  Payne \cite{Payne:06} has shown that $\PP[\Sigma]$ is isomorphic to
  the equivariant Chow ring of~$X_\Sigma$, for \emph{any} fan~$\Sigma$.
  Further developments are documented in \cite{KatzPayne}.
  Proposition~\ref{iso-HT-piecewise-poly} is also valid for more general
  fans, but we have been unable to locate a reference identifying
  $X_\Sigma$ as a $T$-CW~complex; so we proceed by assuming appropriate
  cohomological finiteness conditions, following \cite[Remark
  1.4]{FranzPuppe:exact}.
 \end{remark}

%
%
%
%
%
%
%
%
%

\section{Generators of the ring of piecewise polynomials}\label{generators}

For~$i=0$, \ldots,~$n$ we will write $\sigma_i\in\Sigma$ for the 
full-dimensional cone spanned by all rays except $v_i$.  This cone is
simplicial and full-dimensional because the set~$\{v_j:j\ne i\}$ is a
basis of the $\Q$-vector space $N\otimes_\Z\Q$ for each $i$.
Moreover, given a piecewise polynomial $f$, we will denote the unique
polynomial which coincides with~$f$ on $\sigma_i$ by $f^{(i)}$. We
call a piecewise polynomial \emph{reduced} if it is not divisible in~%
$\PP[\Sigma]$ by any rational prime.

Let $b_{i j}$, where $i\ne j$, be the reduced linear function that
assumes a positive value on $v_i$ and vanishes on all $v_k$, for $i\ne
k\ne j$.

\begin{lemma}\label{value-bij}
  We have that
  \begin{equation}
    b_{i j}(v_i) =\frac{\chi_j}{\gcd(\chi_i,\chi_j)}
    =\frac{\lcm(\chi_i,\chi_j)}{\chi_i}.
  \end{equation}
\end{lemma}
\begin{proof}
  Applying $b_{ij}$ to the relation~(\ref{relation-fan-generators}) yields
  \begin{equation}
    \chi_i\, b_{i j}(v_i) = -\chi_j\, b_{i j}(v_j).
  \end{equation}
  Since $b_{i j}$ is reduced and $v_i$~and~$v_j$ span $N/\ker
  b_{ij}\cong\Z$, the values $b_{i j}(v_i)$ and $b_{i j}(v_j)$ must be
  coprime.  This implies the claimed formula.
\end{proof}

\begin{proposition}\label{generators-linear}
The $b_{ij}$ generate $N^\vee$, the lattice dual to $N$.
\end{proposition}
\begin{proof}
  Multiplying all weights by a constant factor changes neither the
  fan~$\Sigma$ nor the functions~$b_{i j}$, so we may assume that the
  greatest common divisor of the weights equals $1$.  For given~$j$,
  let $N_j$ be the span of the linearly independent
  set~$V_j=\{v_i:i\ne j\}$ and $N_j^\vee$ its dual.  By
  Lemma~\ref{value-bij}, the restriction of each~$b_{ij}$ with~$i\neq
  j$ to~$N_j$ is divisible by a divisor of $\chi_j$, and the quotient
  is an element of the basis dual to~$V_j$.

  Let $M_j<N^\vee$ denote the sublattice generated by those 
  $b_{ij}$ for which $i\neq j$.  Our goal is then to show that
  $M=N^\vee$, where $M$ is generated by the $M_j$.

  We have that 
  \begin{equation}
    N_j^\vee/N^\vee = (N_j^\vee\!/M_j) \bigm/ (N^\vee\!/M_j).
  \end{equation}
  Therefore, the order of $N^\vee\!/M_j$ divides that of
  $N_j^\vee\!/M_j$, which itself divides $\chi_j^n$ from above.
  Hence, the order of $N^\vee\!/M_j$ also divides $\chi_j^n$, and the
  same applies to $N^\vee\!/M$ because
  \begin{equation}
    N^\vee\!/M = (N^\vee\!/M_j) \bigm/ (M/M_j).
  \end{equation}
  This implies that the order of $N^\vee\!/M$ divides the greatest
  common divisor of all~$\chi_j^n$, which we assumed to be $1$.
\end{proof}

Let $a_i$ denote the \emph{Courant function} corresponding
to $v_i$, for $0 \leq i\leq n$.
By this we mean the reduced piecewise linear function that assumes
a positive value on $v_i$ and vanishes on all $v_j$ for $j\ne i$.
Each $\sigma_j$ is simplicial, so $a_i$ is well-defined.

\begin{lemma}\label{generators-piecewise-linear}
  Together with the linear functions, each $a_i$ generates 
  the piecewise linear functions 
  in $\PP[\Sigma]$.
\end{lemma}
\begin{proof}
  Let $f$ be piecewise linear. Then $f-f^{(i)}$ vanishes on~$\sigma_i$,
  and is therefore a multiple of~$a_i$.
\end{proof}

\begin{lemma}\label{aij}\label{aivi}
  We have that
  \begin{equation}
    a_i(v_i)=\frac{\lcm(\chi_0,\dots,\chi_n)}{\chi_i}
    \qquad{\rm and}\qquad
    a_i^{(j)}=\frac{\lcm(\chi_0,\dots,\chi_n)}{\lcm(\chi_i,\chi_j)}\,b_{i j}
  \end{equation}
  in $\PP[\Sigma]$, for all $i\ne j$.
\end{lemma}
\begin{proof}
  By Lemma~\ref{value-bij} we can define a piecewise linear
  function $f$ on $\Sigma$ by setting $f^{(i)}=0$, and $f^{(j)}$ equal
  to the given formula for $j\ne i$. Now let $\chi_k$ be a weight
  with maximal $p$-content, for some prime $p$. If $k=i$, then $p$
  cannot divide $f^{(j)}$ for any~$j\ne k$; and if $k\ne i$, then $p$
  cannot divide $f^{(k)}$. So $f$ is reduced, and $f=a_i$.
\end{proof}

\begin{lemma}\label{aiaj-linear}
  In $\PP[\Sigma]$, the relation
  \begin{equation}\label{relation-bij}
    b_{i j} = \frac{\lcm(\chi_i,\chi_j)}{\lcm(\chi_0,\dots,\chi_n)}\,(a_i-a_j)
  \end{equation}
  holds for all $i\ne j$. 
\end{lemma}

\begin{proof}
  By Lemma~\ref{aij}, we have that
  \begin{equation}
    a_i^{(j)}=-a_j^{(i)}
    = \frac{\lcm(\chi_0,\dots,\chi_n)}{\lcm(\chi_i,\chi_j)}\,b_{i j},
  \end{equation}
  and $a_i^{(i)}=-a_j^{(j)}=0$; so
  $(a_i-a_j)^{(i)}=a_i^{(j)}=(a_i-a_j)^{(j)}$.  Furthermore, by
  Lemma~\ref{generators-piecewise-linear}, $a_i-a_j$ may be written as
  $ma_k+r$ for some integer $m$ and linear function $r$. But every
  $a_k$ restricts to distinct linear functions on maximal cones 
  $\sigma_i$ and $\sigma_j$, so $m=0$. Hence $a_i-a_j$ is linear, 
  and divisible as claimed.
\end{proof}

We now consider higher-degree analogues of the Courant functions~$a_i$.

\begin{lemma}\label{divisibility_a_I}
  For any nonempty subset~$I\subset\{0,\ldots,n\}$, the function
  $\prod_{i\in I} a_i$ is divisible by
  \begin{equation}\label{prodloverl}
    \prod_{i\in I}
      \frac{\lcm(\chi_0,\dots,\chi_n)}{\lcm\{\chi_i,\chi_j:j\notin I\}}
  \end{equation}
in $\PP[\Sigma]$.
\end{lemma}

\begin{proof}
  We look at each prime~$p$ separately.  If the maximal $p$-content
  occurs in $\chi_j$ for some $j\notin I$, then there is nothing to
  prove because the $p$-content of (\ref{prodloverl}) is $1$. We can
  therefore assume that it occurs in $\chi_k$ for some $k\in I$.

  Choose an~$i\in I$ and denote the $p$-contents of $\chi_i$ and
  $\chi_k$ by $q_i$ and $q_k$ respectively.  Then all~$a_i^{(j)}$
  with $j\notin I$ are divisible by~$q_k/q_i$, which is greater than
  or equal to the $p$-content of
  \begin{equation}
    \frac{\lcm(\chi_0,\dots,\chi_n)}{\lcm\{\chi_i,\chi_j:j\notin I\}}.
  \end{equation}
  Taking the product over all~$i\in I$ finishes the proof.
\end{proof}

Hence, for~$I\subset\{0,\ldots,n\}$ we may define the piecewise 
polynomial
\begin{equation}\label{relation-aI}
  a_I = \Bigl(\,\prod_{i\in I}\frac{\lcm(\chi_0,\dots,\chi_n)}
  {\lcm\{\chi_i,\chi_j:j\notin I\}}\Bigr)^{-1}\prod_{i\in I} a_i
\end{equation}
in the $2|I|$-dimensional component of $\PP[\Sigma]$.

\begin{theorem}\label{wproj-generators}
  The ring $H_T^*(\bP(\chi))$ is generated by the functions $a_I$ and
  $b_{ij}$, where $1\le|I|\le n$ and $i\ne j$ respectively. The only
  relations are (\ref{relation-product-all}), (\ref{relation-bij})
  and (\ref{relation-aI}).
\end{theorem}

\begin{proof}
  From (\ref{wpsratcoho}), $H_T^*(\bP(\chi);\Q)$ is generated by the
  $a_i$ subject only to the relation (\ref{relation-product-all}).
  The relations (\ref{relation-bij}) and (\ref{relation-aI}) show that
  the $a_I$ with $|I|>1$ and the $b_{ij}$ are redundant over $\Q$, so
  adding these generators and relations gives an isomorphic ring.
  Since there are no more relations between the $a_I$ and $b_{ij}$ in
  $H_T^*(\bP(\chi);\Q)$, the same is true in $H_T^*(\bP(\chi))$; for
  the latter is free over $\Z$, and injects into
  $H_T^*(\bP(\chi);\Q)$. It remains to show that these elements are
  ring generators.

  By Proposition~\ref{generators-linear}, the $b_{ij}$ generate the
  linear functions, which are the image of~$H^2(BT)$ in
  $H_T^2(\bP(\chi))$.  Hence, by Lemma~\ref{HT-free-generators}, it
  suffices to show that the subgroup generated by the $a_I$ surjects
  onto~$H^*(\bP(\chi))$. In other words, we have to show that $c_m$
  lies in the span of~$\{\iota^*(a_I):|I|=m\}$ for each~$1\le m\le n$.

  For~$m=1$, this is true by Lemma~\ref{generators-piecewise-linear}
  because we know $\iota^*$ itself to be surjective.  Moreover,
  Lemma~\ref{aiaj-linear} implies that all elements $a_i$ are mapped
  to the same element of $H^2(\bP(\chi))$. This must necessarily be a
  generator, which we may assume to be $c_1$.
  
  For~$1<m\le n$, we obtain 
  \begin{subequations}
  \begin{align}
    \iota^*(a_I)
      &= \Bigl(\,\prod_{i\in I}\frac{\lcm(\chi_0,\dots,\chi_n)}{\lcm\{\chi_i,\chi_j:j\notin I\}}\Bigr)^{-1}\prod_{i\in I}\iota^*(a_i) \\
      &= \Bigl(\,\prod_{i\in I}\frac{\lcm(\chi_0,\dots,\chi_n)}{\lcm\{\chi_i,\chi_j:j\notin I\}}\Bigr)^{-1}c_1^m \\
      &= \frac{\prod_{i\in I}\lcm\{\chi_i,\chi_j:j\notin I\}}{\lcm\{\,\prod_{j\in J}\chi_j:|J|=m\,\}}\,c_m
      && \text{by~\eqref{divisibility-powers}}. \label{formula3}
  \end{align}
  \end{subequations}
  We must show that these multiples of~$c_m$ generate
  $H^{2m}(\bP(\chi))$. Once more, we consider each prime $p$
  separately, and let $I$ be the set of indices (which need not be
  unique) that correspond to $m$~weights with greatest $p$-content.
  Since this is also the set~$J$ which maximises the $p$-content of
  the denominator of (\ref{formula3}), we conclude that for each $p$
  there appears a multiple of $c_m$ whose $p$-content is $1$. In other
  words, the greatest common divisor of all multiples is $1$, as
  required.
\end{proof}

%
%
%
%
%
%
%
%
%

\section{An example}\label{example}

We illustrate the results of the preceding section in the case
$\chi=(1,2,3,4)$, which confirms that the elements $b_{ij}$ cannot be
omitted from the statement of Theorem~\ref{wproj-generators}.

We choose $v_0=(-2,-3,-4)$, $v_1=(1,0,0)$, $v_2=(0,1,0)$, and
$v_3=(0,0,1)$.  So isomorphisms $(S^1)^4/S^1\langle
1,2,3,4\rangle\leftrightarrow(S^1)^3$ identifying the torus $T$ of
(\ref{wtdtorus}) are induced by
$(t_0,t_1,t_2,t_3)\mapsto(t_0^{-2}t_1,t_0^{-3}t_2,t_0^{-4}t_3)$ and
$(u_1,u_2,u_3)\mapsto(1,u_1,u_2,u_3)$ respectively, and its action on
$\bP(1,2,3,4)$ is equivalent to that of $(S^1)^3$, given by
\begin{equation}\label{act1234}
  (u_1,u_2,u_3)\cdot[x_0,x_1,x_2,x_3]=
  [x_0,u_1x_1,u_2x_2,u_3x_3]
\end{equation}
on homogeneous coordinates.

Writing an element $f$ of $\PP[\Sigma]$ as
$f=(f^{(0)},f^{(1)},f^{(2)},f^{(3)})$, and the canonical basis of
$N^\vee$ as $(x,y,z)$, we deduce
\begin{subequations}
\begin{align}
  a_0 &= (0, -6x, -4y, -3z), \\
  a_1 &= (6x, 0, 6x-4y, 6x-3z), \\
  a_2 &= (4y, -6x+4y, 0, 4y-3z), \\
  a_3 &= (3z, -6x+3z, -4y+3z, 0).
\end{align}
\end{subequations}
Each $a_i$ is reduced, although individual components $a_i^{(j)}$ may
have non-trivial divisors.  The situation changes for products~%
$\prod_{i\in I}a_i$, because components \smash{$a_i^{(j)}$} with~$j\in
I$ are multiplied by $0$; this is the essence of
Lemma~\ref{divisibility_a_I}. Since $\lcm\{1,2,3,4\}=12$, we obtain
\begin{subequations}
\begin{gather}
  a_{01} = a_0 a_1, \quad
  a_{02} = a_0 a_2 / 3, \quad
  a_{03} = a_0 a_3 / 2, \\
  a_{12} = a_1 a_2 / 3, \quad
  a_{13} = a_1 a_3 / 4, \quad
  a_{23} = a_2 a_3 / 6, \\
  a_{012} = a_0 a_1 a_2 / 9, \quad
  a_{013} = a_0 a_1 a_3 / 8, \\
  a_{023} = a_0 a_2 a_3 / 36, \quad
  a_{123} = a_1 a_2 a_3 / 72.
\end{gather}
\end{subequations}

Any globally linear function obtained from the $a_i$ is a linear 
combination of
\begin{subequations}
\begin{gather}
  a_1-a_0=6x, \quad
  a_2-a_0=4y, \quad
  a_3-a_0=3z, \label{sub1}\\
  a_2-a_1=4y-6x, \quad
  a_3-a_1=3z-6x, \quad
  a_3-a_2=3z-4y,\label{sub2}
\end{gather}
\end{subequations}
by Proposition \ref{generators-linear} and Lemma \ref{aiaj-linear}.
The functions $b_{10}=x$, $b_{20}=y$ and $b_{30}=z$ are obviously not
in the span of the $a_i$. On the other hand, given $x$, $y$ and $z$,
the remaining $b_{ij}$ are redundant, as are most of the $a_I$. In
fact we may write
\begin{equation}\label{coh1234}
H_T^*(\bP(1,2,3,4))\cong\Z
\big[x,y,z,a_3,a_{23},a_{123}\big]\big/(a_0a_1a_2a_3)
\end{equation}
by Lemma \ref{HT-free-generators}, subject to the relations above. But
we know of no canonical choice for a minimal set of generators.

For examples in which $\chi_i\neq 1$ for any $i$, the construction of
a fan $\Sigma$ and an explicit $T$-action on $\bP(\chi)$ is more
difficult, and amounts to completing $\chi/\gcd(\chi_0,\ldots,\chi_n)$
to a lattice basis.  Further details may be found in
\cite[\S4.2]{Kasprzyk:06}.

%
%
%
%
%
%
%
%
%

\section{Recovering the weights}\label{recover-weights}

It is clear from the definition
(\ref{definition-weighted-projective-space}) that $\bP(\chi)$ does not
change if all weights are multiplied by the same factor.  In
particular, we may always divide the weights by their greatest common
divisor, as in the proof of Proposition~\ref{generators-linear}.
Moreover, if every weight except $\chi_i$ is divisible by some prime
$p$, then $\bP(\chi)$ is equivariantly isomorphic to $\bP(\chi')$,
where $\chi'=(\chi_0/p,\dots,\chi_{i-1}/p,\,\chi_{i},\,\chi_{i+1}/p,
\dots,\chi_n/p)$ \cite{Dolgachev:82}. This may be seen from the toric
viewpoint: for (\ref{relation-fan-generators}) implies that $v_i$ is
also divisible by~$p$, and continues to hold when $v_i$ and $\chi$ are
replaced by $v_i/p$ and $\chi'$ respectively. So the fan $\Sigma$ is
unchanged, and the corresponding toric varieties are isomorphic.

By repeating these simplifications, we can always ensure that two or
more weights are not divisible by $p$, for each prime $p$.  The
resulting weight vector is uniquely defined by~$\chi$ (up to order),
and the corresponding weights are called \emph{normalised}.

\begin{theorem}\label{cohodetermineschi}
The graded ring~$H_T^*(\bP(\chi))$ determines the normalised weights.
\end{theorem}
\begin{proof}
  The length $n+1$ of $\chi$ is given by the rank of the free abelian
  group $H_T^2(\bP(\chi))$.  So we may interpret $H_T^*(\bP(\chi))$ as
  a ring of piecewise polynomials $\PP[\Sigma]$, where $\Sigma$ has
  cones $\sigma_i$ and Courant functions $a_i$, for $0\leq i\leq n$.

  According to relation (\ref{relation-product-all}), we may choose
  piecewise linear functions $f_0$, \dots, $f_n$ in $H_T^2(\bP(\chi))$
  that are reduced, non-zero, and satisfy $f_0\dots f_n=0$. On each
  cone~$\sigma_j$, some $f_i$ must therefore vanish; but it cannot
  vanish on $\sigma_k$ for any other $k\ne j$, or else it would also
  vanish on every ray of $\Sigma$ and be identically zero. Because
  $f_i$ is reduced, it follows that $f_i=\pm a_j$. So we may assume
  that $f_i=\pm a_i$ for every~$0\leq i\leq n$, by permuting the cones
  as necessary.
  
  Given any $i\neq j$, we may now read off the $p$-content~$q_{ij}$
  of~$a_i-a_j$ from the functions $f_i$, as follows. Since $f_i=\pm
  a_i$ and $f_j=\pm a_j$, we know that $q_{ij}$ is the $p$-content of
  either~$f_i-f_j$ or~$f_i+f_j$; in fact it is the larger of the two
  (and the smaller is the $p$-content of $a_i+a_j$). This is
  because $a_j$ restricts to $0$ on $\sigma_j$, whence
  \smash{$(a_i-a_j)^{(j)}=(a_i+a_j)^{(j)}$}. But $a_i-a_j$ is globally
  linear by Lemma \ref{aiaj-linear}, so its $p$-content is unaltered
  by restriction to $\sigma_j$, whereas that of $a_i+a_j$ may
  increase.

  Appealing to Lemma~\ref{aiaj-linear} once more, we find that
  $\lcm(\chi_0,\dots,\chi_n)/\lcm(\chi_i,\chi_j)$ has $p$-content
  $q_{ij}$.  Moreover, there exist integers $j$ and $k$ for which
  $\lcm(\chi_j,\chi_k)$ is not divisible by~$p$, since the weights are
  normalised.  So for $i\neq j$ or $k$, the $p$-content of~$\chi_i$ is
  precisely $q_{jk}/q_{ik}$, and the weights are completely
  determined up to order.
\end{proof}

\begin{remark}
  The analogue of Theorem \ref{cohodetermineschi} is false for
  ordinary cohomology, because the divisibility
  rule~(\ref{divisibility-powers}) does not take into account the
  relationship between the distribution of the $p$-contents of the
  weights for different primes $p$. For example, the graded rings
  $H^*(\bP(1,2,3))$ and $H^*(\bP(1,1,6))$ are isomorphic, since
  $c_1^2=6c_2$ in both cases. However, $\bP(1,2,3)$ and $\bP(1,1,6)$
  cannot be homeomorphic, because the former has two singular points,
  and the latter only one; nevertheless, they are both homotopy
  equivalent to the $2$-cell complex $S^2\cup_{6\eta}e^4$, where
  $\eta$ generates $\pi_3(S^2)\cong\Z$.

  From the toric viewpoint, both quotient spaces $\bP(1,2,3)/T$ and
  $\bP(1,1,6)/T$ may be identified with the $2$-simplex. On the other
  hand, it follows from Theorem \ref{cohodetermineschi} that the
  corresponding homotopy quotients cannot even be homotopy equivalent.
\end{remark}

%
%
%
%
%
%
%
%
%

\section{Weighted projective bundles}\label{projective-bundles}

Suppose given complex line bundles $L_{i}$ over a base space $X$ for
$0\leq i\leq n$, and denote their direct sum by
$D=L_0\oplus\dots\oplus L_n$. The torus $T'=(S^1)^{n+1}$ acts on~$D$,
and on the corresponding sphere bundle~$S(D)$, in canonical
fashion. The associated \emph{weighted projective bundle} over $X$ has
fibre $\bP(\chi)$, and total space the quotient
\begin{equation}
  \bP(D,\chi)=S(D)/S^1\langle\chi_0,\ldots,\chi_n\rangle.
\end{equation}
The universal example is given by $E=ET'\times_{T'}\C^{n+1}$ over
$BT'$. The reasoning of Lemma~\ref{HT-free-generators} shows that
$H^*(\bP(E,\chi))$ is a free $H^*(BT')$-module of rank $n+1$, and the
naturality of the Serre spectral sequence implies the same result for
arbitrary~$D$.

If all weights are equal to~$1$, then $\bP(D,\chi)$ is an ordinary
projective bundle, and $H^*(\bP(D,\chi))$ is generated by the Chern
classes $c_1(L_i)$ and $-\xi$, subject only to the relation
(\ref{chern-projective-prod}). In \cite[Ch.III]{AlAmrani:97}, 
Al~Amrani generalised (\ref{chern-projective-prod}) to those weighted 
projective bundles whose $\chi_i$ form a divisor chain. In all such 
cases, he proved that
\begin{equation}\label{chern-weighted-projective}
  \prod_{i=0}^n\Bigl(\xi
    +\frac{\lcm(\chi_0,\dots,\chi_n)}{\chi_i}\,c_1(L_i)\Bigr) = 0
\end{equation}
for a certain $\xi\in H^2(\bP(D,\chi))$, which restricts to~$c_1\in
H^2(\bP(\chi))$ on fibres. By naturality, it suffices to verify
(\ref{chern-weighted-projective}) for the universal case; in other
words, in $H_{T'}^*(\bP(\chi))$, where $T'$ acts on $\bP(\chi)$ via
the projection
  \begin{equation}\label{T-projection}
    T'\longrightarrow T= 
    T'/S^1\langle\chi_0,\ldots,\chi_n\rangle.
  \end{equation}

We may profitably describe $H_{T'}^*(\bP(\chi))$ in terms of  
piecewise polynomials, as follows. 

Let $\pi\colon N'=\Z^{n+1}\to N$ be the epimorphism defined on
canonical basis vectors by $\pi(e_i)=v_i$, for $0\leq i\leq n$. Its
kernel is the subgroup generated by $\chi_0 e_0+\dots+\chi_n e_n$,
which we abbreviate to $u$. The pull-back
$\Sigma'=\{\pi^{-1}(\sigma):\sigma\in\Sigma\}$ of $\Sigma$ is a
\emph{generalised fan} \cite{Brion:96} because every cone contains
the line through $u$, and $\pi$ induces a monomorphism
$\pi^*\colon\PP[\Sigma]\to\PP[\Sigma']$ of piecewise polynomial rings.
\begin{lemma}
As $H^*(BT')$-algebras, $H_{T'}^*(\bP(\chi))$ is naturally 
isomorphic to $\PP[\Sigma']$; furthermore, the projection 
$T'\to T$ induces the homomorphism
$\pi^*\colon H_T^*(\bP(\chi))\to H_{T'}^*(\bP(\chi))$. 
\end{lemma}
\begin{proof}
  On cohomology rings we have the natural monomorphism
  \begin{equation}
    H_T^*(\bP(\chi))\longrightarrow H_{T'}^*(\bP(\chi))
      =H_T^*(\bP(\chi))\otimes H^*(BS^1),
  \end{equation}
  because $T'$ splits as $T\times S^1$, where
  $S^1$ is the kernel of (\ref{T-projection}) and acts 
  trivially on~$\bP(\chi)$. So the freeness of~$H_T^*(\bP(\chi))$
  over~$H^*(BT)$ implies the freeness of $H_{T'}^*(\bP(\chi))$ 
  over~$H^*(BT')$. Moreover, every isotropy subgroup of the
  $T'$-action is connected.
  
  We can therefore imitate the reasoning of
  Proposition~\ref{iso-HT-piecewise-poly}. For cones
  ~$\sigma'\in\Sigma'$ that have $k=\codim\sigma'\le1$ (and in fact
  for all cones ~$\sigma'$), we have that
  \begin{equation}
    H_{T'}^{*+k} 
      (\bar\cO_{\sigma'},\partial\cO_{\sigma'})
      \cong \Z[\sigma'], 
  \end{equation}
  and the boundary map in cohomology corresponds to the 
  restriction of polynomials up to sign. Hence the integral 
  Chang--Skjelbred sequence is exact, which in turn identifies 
  $H_{T'}^*(\bP(\chi))$ with $\PP[\Sigma']$. 

  Since for all~$\sigma\in\Sigma$ 
  and~$\sigma'=\pi^{-1}(\sigma)$ the map 
  \begin{equation} 
    H_T^*(\bar\cO_\sigma,\partial\cO_\sigma)\to 
    H_{T'}^*(\bar\cO_{\sigma'},\partial\cO_{\sigma'})
  \end{equation} 
  corresponds to the pull-back of functions by $\pi$,
  the same applies to the restriction of 
  equivariant cohomology from $T$ to~$T'$. 
\end{proof}

Because no cone of $\Sigma'$ contains every $e_i$, we can define $\xi$
as the piecewise linear function on $N'$ which takes the values
$\xi(u)=-\lcm(\chi_0,\dots,\chi_n)$, and $\xi(e_i)=0$ for all $0\leq
i\leq n$. Equivalently,
\begin{equation}
  \xi^{(i)}=-\frac{\lcm(\chi_0,\dots,\chi_n)}{\chi_i}\,x_i
  \qquad
\end{equation}
for all $i$, where $(x_i)$ denotes the basis dual to 
$(e_i)$ for $N'$. 

\begin{theorem}\label{weighted-chern} 
  As an element of $H_{T'}^*(\bP(\chi))$, the cohomology 
  class $\xi$ restricts to~$c_1$ in $H^2(\bP(\chi))$ 
  and satisfies equation~(\ref{chern-weighted-projective}). 
\end{theorem}
\begin{proof}
  As elements of $\PP[\Sigma']$, we have that 
  \begin{subequations}\label{pull-back-ai}
  \begin{align}
      \pi^*(a_i)^{(j)} &= \frac{\lcm(\chi_0,\dots,\chi_n)}{\chi_i}\,x_i
                         -\frac{\lcm(\chi_0,\dots,\chi_n)}{\chi_j}\,x_j \\
                       &= \frac{\lcm(\chi_0,\dots,\chi_n)}{\chi_i}\,x_i
                         +\xi^{(j)},
  \end{align}
  \end{subequations}
  because the right-hand side vanishes on~$u$ and assumes the value
  $a_i(v_k)$ on $e_k$, for all~$k$.  Since $\xi$ differs
  from~$\pi^*(a_i)$ by a linear function, it restricts to the same
  element in~$H^2(\bP(\chi))$ as $\pi^*(a_i)$~and~$a_i$, namely
  $c_1$.  Identifying $c_1(L_i)$ with $x_i$, we conclude that
  equation~(\ref{chern-weighted-projective}) is nothing but the
  pull-back of relation~(\ref{relation-product-all}).
\end{proof}
No other relation is required to describe $H^*(\bP(E,\chi))$
rationally, nor therefore integrally.

%
%
%
%
%
%
%
%
%

\bigskip
\footnotesize

\noindent
\textsc{Department of Mathematics, Rider University,
  Lawrenceville NJ, 08648, U.S.A.}

\noindent
\emph{E-mail address:}~\texttt{bahri@rider.edu}

\smallskip

\noindent
\textsc{Fachbereich Mathematik, 
  Universit\"at Konstanz, 78457 Konstanz, Germany}

\noindent
\emph{E-mail address:}~\texttt{matthias.franz@ujf-grenoble.fr}

\smallskip

\noindent
\textsc{School of Mathematics, University of Manchester,
  Oxford Road, Manchester M13 9PL, England}

\noindent
\emph{E-mail address:}~\texttt{nige@maths.manchester.ac.uk}

\end{document}